\documentclass[12pt]{amsart}
\usepackage{amsmath}
\usepackage{amssymb}
\usepackage{amsthm}
\usepackage{graphicx}
\usepackage{verbatim}
\usepackage{enumerate}
\usepackage[margin=1in]{geometry}
\usepackage{listings}
\usepackage{todonotes}
\theoremstyle{plain}
\newtheorem{thm}{Theorem}[section]

\newtheorem{lem}[thm]{Lemma}

\newtheorem{defn}[thm]{Definition}

\theoremstyle{remark}

\newtheorem{remark}{Remark}

\newcommand{\pois}{\mathcal{P}}

\def\Z{\mathbb{Z}}
\def\Q{{\bf Q}}
\def\R{\mathbb{R}}

\def\S{{\mathcal S}}

\def\ee{\varepsilon}
\def\E{{\mathbb E}}
\def\P{{\mathbb P}}
\def\Cox{\hfill \Box}

\def\disp{\displaystyle}
\def\one{{\bf 1}}
\def\|{{\, | \, }}
\def\F{{\mathcal F}}

\def\X{{\bf X}}
\def\Y{{\bf Y}}
\def\yy{{\bf y}}

\def\invt{{\mathcal I}}

\def\pt{{\tilde{p}}}
\def\qt{{\tilde{q}}}

\def\rate{\eta}
\def\tphi{\phi_{W,n}}
\def\tpsi{\psi_{W,n}}

\DeclareMathOperator{\sumset}{sumset}
\begin{document}
\title[Four permutations]{Four random permutations conjugated by an adversary
generate $\S_n$ with high probability}

\author{Robin Pemantle}
\address{Department of Mathematics,
University of Pennsylvania,
209 South 33rd Street,
Philadelphia, PA 19104, USA}
\email{pemantle@math.upenn.edu}
\author{Yuval Peres}
\address{Microsoft Research,
1 Microsoft Way, Redmond, WA, 98052, USA}
\email{peres@microsoft.com}
\author{Igor Rivin}
\address{Temple University,
 1805 N Broad St, Philadelphia, PA}
\curraddr{Mathematics Department, Brown University}
\email{igor.rivin@temple.edu}
\begin{abstract}
We prove a conjecture dating back to a 1978 paper of
D.R.\ Musser~\cite{musserirred}, namely that four random
permutations in the symmetric group $S_n$ generate a transitive
subgroup with probability $p_n > \ee$ for some $\ee > 0$ independent
of $n$, even when an adversary is allowed to conjugate each of the
four by a possibly different element of $\S_n$ (in other words,
the cycle types already guarantee generation of $\S_n$).  This is
closely related to the following random set model.  A random set
$M \subseteq \Z^+$ is generated by including each $n \geq 1$
independently with probability $1/n$.  The sumset $\sumset(M)$ is formed.
Then at most four independent copies of $\sumset(M)$ are needed before
their mutual intersection is no longer infinite.
\end{abstract}

\thanks{Igor Rivin would like to thank the Brown University Mathematics Department and ICERM for their hospitality and financial support during the preparation of this paper. Robin Pemantle was supported in part by NSF grant  \# DMS-1209117}
\maketitle

\setcounter{equation}{0}
\section{Introduction}
\label{sec:intro}

\subsection{Background and motivation}
The roots of this work are in computational algebra.  It is
a result going back to van der Waerden~\cite{van1934seltenheit}
that most polynomials $p(x) \in \Z [x]$ of degree $n$ have
Galois group $S_n$.  Computing the Galois group is a
central problem in computational number theory and is a
fundamental building block for the solution of seemingly
unrelated problems (see~\cite{rivin2013large} for an extensive
discussion).  Therefore, one cannot take for granted being in
the ``generic'' case and one would like an effective and
speedy algorithm for determining whether the Galois group of $p(x)$
is the full symmetric group.

There are deterministic polynomial time algorithms to answer
this.  The first is due to S.\ Landau;  a simpler and more 
efficient algorithm was proposed by the third author 
(see~\cite{rivin2013large}).  These algorithms, however, are of 
purely theoretical interest due to their very long run times 
(their complexity is of the order of $O(n^{40}),$ where $n$ 
is the degree of the polynomial).
The best algorithms in practice are Monte Carlo algorithms.  To discuss
Monte Carlo testing for full Galois group, one begins with two classical
results\footnote{In the literature, the much harder Chebotarev Density
Theorem is often used in place of the Frobenius Density Theorem.}.

\begin{thm}[Dedekind]
\label{dedthm}
If $p(x)$ is square-free modulo a prime $q$ and the factorization of
$p(x)$ modulo $q$ into irreducible factors yields degrees
$d_1, d_2, \ldots, d_k$, then the Galois group of $G$ $p(x)$
has an element whose cycle decomposition has lengths precisely
$\{ d_1, \ldots, d_k \}$.
\end{thm}

\begin{thm}[Frobenius Density Theorem]
\label{frobthm}
The density of prime numbers $q$ for which $p(x)$ mod $q$ has factors
whose degrees are $d_1 , \ldots , d_k$ is equal to the density in
the Galois group $G \subseteq \S_n$ for $p(x)$ of elements of $\S_n$
with cycle type $d_1, \ldots, d_k$.
\end{thm}
\begin{remark}
Theorem \ref{frobthm} is useless without effective convergence bounds. The first step in this direction was made by the J. Lagarias and A. Odlyzko \cite{lagariasodlyzko} -- they proved \emph{conditional} (on the Riemann hypothesis for certain L-functions) results with ``effectively computable'' (but quite hard to compute) constants. A couple of years later, Oesterl\'e \cite{oesterle1979versions} claimed a computation of the constants, but his computation has not been published in the intervening 35 years (despite being used by J.-P.~Serre in \cite{serre1981IHES}). Finally, the problem was put to rest by B.~Winckler in \cite{winckler2013th}) at the end of 2013(!) -- Winckler shows both unconditional and conditional results (with somewhat worse constants  in the latter case than those claimed by Oesterl\'e).
\end{remark}
Together, these two results tell us that without yet knowing $G$ we can
uniformly sample cycle decompositions $V_i = \{ d_{i,1} , \ldots ,
d_{i,k(i)} \}$ of elements of $G$ by sampling integers $q_i$ at
random and setting $V_i$ equal to the set of degrees of the irreducible
factors of $p(x)$ modulo $q_i$ (it should be noted that factoring modulo a prime can be done quite efficiently using variants of Berlekamp's algorithm).  A result of C.\ Jordan allows us to
turn this into a probabilistic test for $G = \S_n$ with certain
acceptance and possible false rejection.

\begin{thm}[C. Jordan]
\label{jord2thm}
Suppose a subgroup $H$ of $\S_n$ ($n > 12$) acts transitively
on $[n]$.  If it contains at least one cycle of prime length
between $n/2 + 1$ and $n-5$, then it is either $\S_n$ or
the alternating group $\mathcal{A}_n.$
\end{thm}

Certification that $G$ is not alternating and contains
at least one long prime cycle is trivial: we just check
that at least one of the lists $V_1 , \ldots V_r$ corresponds
to an odd permutation class and at least one contains a
prime value in $[n/2+1 , n-5]$ -- some power of the corresponding
permutation will be a long prime cycle.  In a uniform random
permutation, the cycle containing a given element, say~1, has
length exactly uniform on $[n]$.  The Prime Number Theorem
guarantees that the number of primes in $[n/2+1 , n-5]$ is
asymptotic to $n / (2 \log n)$.  It follows that if $G$ is truly $\S_n$,
then each $V_i$ contains a large prime with probability at least
$(1 + o(1)) / (2 \log n)$.  Also, each $V_i$ corresponds to an odd
class with probability $1/2$.  Therefore, if $G$ is truly $\S_n$,
we will quickly discover that the hypotheses of Theorem~\ref{jord2thm}
other than transitivity are satisfied.

Establishing transitivity of $G$ when we know only $V_1 , \ldots , V_r$
must involve showing that {\em any} set of cycles in these
respective conjugacy classes generates a transitive subgroup
of $\S_n$.  Let us say in this case that classes $V_1 , \ldots ,
V_r$ {\bf invariably generate} a transitive group.
If the action of $G$ leaves a subset $I$ of $[n]$
invariant, then $|I|$ will appear as a sum of cycle sizes of
every element of $G$.  The converse holds as well: if the
sumsets of $V_1 , \ldots , V_r$ have no common intersection
then the corresponding permutations invariably generate 
a transitive group.
This leads to the following test:

\begin{quote}
{\bf Algorithm:} Sample some random primes $\{ q_1, q_2, \ldots q_r \}$,
compute the degree sets $V_i := \{ d_{i,1} , \ldots , d_{i,k(i)} \}$
of the factors of $p$ modulo $q_i$, and the sumsets $S_i := \sumset(V_i)$.
If the sets $S_i$ have some element in common other than~0 and $n$,
or if none of the sets $V_i$ contains a prime grater than $n/2$,
or if all $r$ conjugacy classes are even, then output NEGATIVE,
otherwise output POSITIVE.
\end{quote}

In the algorithm above, we have implicitly defined \emph{sumset}: 
\begin{defn}
The \emph{sumset} of a (multi)set $S=\{k_1, \ldots, k_l\}$ is the (multi)set of all sums of subsets of $S.$
\end{defn}

What needs to be checked next is that that we can choose $r (\ee)$
not too large so that if $p(x)$ does have full Galois group then
a NEGATIVE output has probability less than $\ee$.  For this we
need to answer the question: given $\ee > 0$, how many uniformly
random permutations in $\S_n$ do we have to choose before their
cycle length sumsets have no common value in $\{ 1 , \ldots , n-1 \}$?
If there is a number $m_0$ such that this probability is at least
$\delta$ for $m_0$ random permutations, then it is at least
$1 - (1-\delta)^j$ for $jm_0$ permutations.  Therefore we may
begin by asking about the value of $m_0$: how many IID uniform
permutations are needed so that their cycle length sumsets have
no nontrivial common value with probability that remains bounded
away from zero as $n \to \infty$?

\begin{figure}
\label{expres}
\caption{Experimental results}
\includegraphics[width=6.5in]{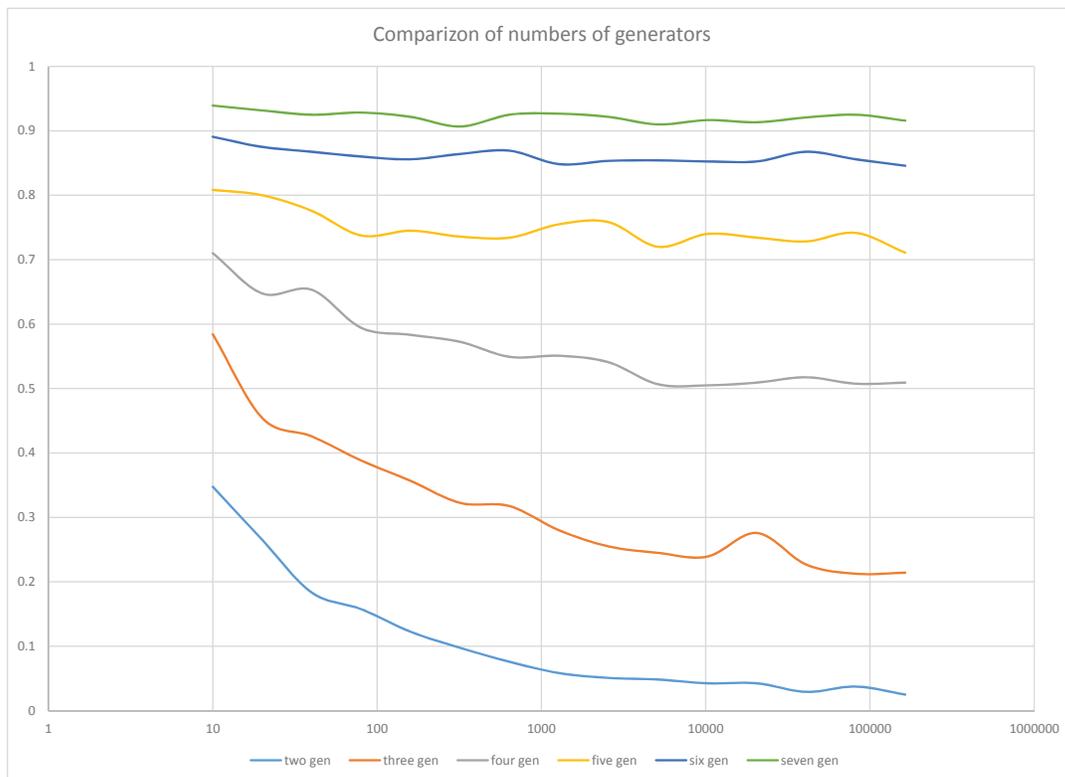}
\end{figure}

It turns out that this question was first raised by
D.~R.~Musser~\cite{musserirred}, for reasons similar to ours.
Musser did some experiments (where he was hindered both by the
performance of the hardware of the time and by using an algorithm
considerably inferior to the one we describe below), and observed
that $5$ elements should be sufficient; see also~\cite{davsmith}.
More modern experimental evidence (See Figure \ref{expres}) is as follows.
Each curve represents the probability that some number of random
elements of $S_n$ invariably generates a transitive subgroup,
where the $x$-axis measures $n$ logarithmically.  The goal is to
prove that one of these curves does not go to zero as $n \to \infty$.
Evidently, even the lowest of these curves does not seem to go
to zero very fast (the horizontal axis is logarithmic), thus we
might believe the question to be delicate.

This question (again, for Galois-theoretic reasons) was considered
by J.~Dixon in his 1992 paper~\cite{dixonrand}, and he succeeded
in showing that $O(\sqrt{\log n})$ elements are sufficient for
fixed $\ee$.  Pictorially, to get above $\ee$ on the graph, it
would suffice to go to the curve numbered $C \sqrt{\log n}$ from
the bottom.  He conjectured, as did Musser, that his bound
was not sharp, and $O(1)$ elements should suffice.  He proved
that if that is, indeed, true, then to check that the Galois group
is all of $S_n$ we need to factor modulo $O(\log \log n)$ primes.
Dixon's $O(1)$ conjectured was proved by T.~Luczak and L.~Pyber in 1993
(\cite{luczakPyber}), however the implied constant was absurdly high:
on the order of $2^{100}$ (and it can be shown that their method cannot
be improved to yield a qualitatively better result).

\subsection{Main results}

Our main result is that $m_0 \leq 4$.
We do not settle whether $m_0$ could be~2 or~3, though we
discuss why very likely $m_0 = 4$ (though experimental evidence
is inconclusive) and why proving this via analyses such as ours 
would require significantly more work.


Let $\P_N$ denote the uniform measure on the symmetric group, $\S_N$.
For a permutation $\sigma \in \S_n$, let $\invt (\sigma)$ denote the set
of sizes of invariant sets of $\sigma$, that is,
$$\invt(\sigma) := \{ |I| : I \mbox{ is a proper subset of } [N]
   \mbox{ and } \sigma [I] = I \} \, .$$
In other words, $\invt(\sigma) = \sumset (V(\sigma))$ when $V$ is
the multiset of cycle lengths of $\sigma$.
Trivially, the set $\invt(\sigma)$ is symmetric about $N/2$, meaning
that it is closed under $k \mapsto N-k$.
As usual, $\P_N^j$ denotes the $j$-fold product of uniform measures
on $\S_N$.

\begin{thm}[Main result] \label{th:main}
There is a positive number $b_0$ such that for all $N$,
$$\P_N^4 \left \{ (\sigma_1 , \sigma_2 , \sigma_3 , \sigma_4) :
   \bigcap_{j=1}^4 \invt(\sigma_j) = \emptyset \right \} \geq b_0 \, .$$
\end{thm}


The ideas behind the proof of this are more evident when we take
$N$ to infinity, resulting in the following Poisson model.
Let $\P$ denote the probability measure on $(\Z^+)^\infty$ making
the coordinates $X_j (\omega) := \omega_j$ into independent
Poisson variables with $\E X_j = 1/j$.
Let $M = M(\omega)$ be the multiset having $X_k$ copies of the
positive integer $k$.  Let $S = S(\omega) = \sumset (M(\omega))$
be the sumset; we may define this formally by
$$S = \left \{ \sum_k a_k \cdot k : a_k \leq X_k \mbox{ for all } k
   \right \} \, .$$
This is the analogue in the Poisson model of the set $\invt(\sigma)$
of sums of cycle lengths in the group theoretic model.

Let $\P^4$ denote the fourfold product of $\P$ on $((\Z^+)^\infty)^4$
and for a 4-sequence $(\omega^1, \omega^2, \omega^3 \omega^4) \in
((\Z^+)^\infty)^4$, let $X_{r,k}$ denote the $k^{th}$ coordinate of
$\omega^r$.  Let $S(\omega_r)$ denote the set of sumsets of
$\omega^r$.  Our main result on the Poisson model is:

\begin{thm}[Poisson result] \label{th:poisson}
$$\P^4 \left ( \bigcap_{r=1}^4 S(\omega^r) = \emptyset \right ) > 0 \, .$$
\end{thm}

We require a number of estimates of probabilities associated with
the random sumset $S$.  The most straightforward quantity to define,
though, as it turns out, not the most useful, is the marginal probability
$p_n  := \P (n \in S)$ of finding a number $n$ in the random sumset.
This quantity is estimated as follows.
\begin{thm}[marginal probabilites] \label{th:p_k}
Let $\disp{\rate = \frac{1 - \log 2 - \log(1/\log 2)}{\log 2}}
\approx -0.08607\ldots$. \hfill \\
Then $p_n = n^{\rate + o(1)}$.
\end{thm}

We remark that the exponent $\rate$ is familiar from number theoretic
contexts.  For example, the asymptotic density of integers $m$ having
a divisor in the interval $[N,2N]$, is a quantity $g(N)$ known to
satisfy $g(N) \sim (\log N)^{\rate + o(1)}$ as $N \to \infty$
(see, e.g.,~\cite{HaTe1988}).

\subsection{Discussion}

The analysis relies on the following lemma of Arratia and Tavar{\'e},
to the effect that the joint distribution of number of cycles of
lengths up to $m = o(N)$ of a random permutation of $\S_N$ look
like independent Poissons (see also~\cite{Gran2006} for further 
refinements).

\begin{lem}[\protect{\cite[Theorem~2]{arratia-tavare}}] \label{lem:AT}
Let $Q_{N,m}$ be the joint distribution, for $1 \leq k \leq m$,
of the number of $k$-cycles in a uniform random permutation
in $\S_N$.  Let $\nu_m := \prod_{j=1}^m \pois (1/j)$ denote the
product of Poisson laws with respective means $1/j$.  Then
there is a constant $C>0$ such that the \textbf{total variation distance}
between these two distributions is bounded above by
$$||Q_{N,m} - \nu_m||_{TV} \leq \exp (- C (N/m) \log (N/m)) \, .$$
In particular, $||Q_{N,m} - \nu_m||_{TV} \to 0$ as $N/m \to \infty$.
\end{lem}

Our main result is proved by showing that the random set $\invt(\sigma)$
behaves roughly like a set of dimension $\ln 2$, that is,
it typically has density $n^{\ln 2 - 1 + o(1)}$ near $n$.
It follows that interseting four of these yields a co-dimension
greater than $1$, which is characteristic of a random set which
is almost surely finite and possibly empty.

Given the relatively clean Poisson approximation, one might wonder
why there is any difficulty at all in proving such a result.
The reason for the difficulty is that the averages of certain
quantities are dominated by exceptionally large contributions 
from sets of small probability and therefore do not represent
the typical values.  For example, let $q_{n,k}$ be the probability
that there is an invariant set of size $k$ and let $e_{n,k}$ be the 
expected number of invariant sets of size $k$.  Because $q_{n,k} \ll 1$,
one might expect that $e_{n,k} \approx q_{n,k}$, but it turns out 
that $e_{n,k} = 1$ precisely, for all $n$ and $k$ (simply check that 
each $k$-set has probability $\disp{\binom{n}{k}^{-1}}$ of being an 
invariant set).  Thus $\P_N (k \in \invt(\sigma))$ is much smaller 
than the expected number of representations of $k$ as a sum of cycle 
lengths.  A similar phenomenon holds for the Poisson model.  
The expected number of ways that the integer $n$ is the sum of 
elements of the random multiset $M$ is the $z^n$ coefficient 
in the generating function $\disp{\prod_{k=1}^\infty \exp (z^k / k)}$, 
which simplifies to precisely~1.  We see that $p_n$ is much smaller 
than this expectation.

What  is more subtle is that even $p_n$ does not give the
right estimate.  The right estimate is what is known in the
statistical physics as the {\bf quenched} estimate.  This
is the estimate obtained when a $o(1)$ portion of the probability
space is excluded which contributes non-negligibly to the
quantity in question, in this case $p_n$.  Holding key parameters
at their typical values produces a ``correct'', quenched estimate,
$\pt_n$.  It may sound strange to ask 
what is the probability that $n$ is in the sumset under 
typical behavior because $p_n$ is already a probability,
meaning it is averaged over all behaviors.  To be clear, 
to obtain the quenched estimate $\pt_n$, we exclude a set 
of arbitrarily small probability (but a single set for all $n$), 
such that off of this set the probability $\pt_n$ of finding 
$n \in S$ is much smaller than $n^{-\rate}$, decaying instead 
like $n^{\log 2 - 1}$.

This is important because $|\rate| = 0.08607$ is a bit larger than
$1/12$, whereas $1 - \log 2$ is a little larger than $1/4$.  Showing
that a set has co-dimension $|\rate|$ indicates that one should intersect
twelve independent copies in order to arrive at the empty set.
When the random sets have co-dimension $1 - \log 2$, however, only
four should be required.  Interestingly, it is no easier to prove
that 12 suffice than that 4 suffice, because the estimate of $p_n$
is as hard as the estimate of $\pt_n$.  

Finally, we note that this is in some sense the ``easy'' direction.
To show that the fourfold intersection is finite in the Poisson
model and often empty in the permutation model requires
only an upper bound on the marginals $p_k$.
To show that a threefold intersection does
not suffice would require an upper bound on the probability
of $j$ and $k$ both being in $\invt(\sigma)$.  This appears more
difficult.

\setcounter{equation}{0}
\section{Estimates for the Poisson model}
\label{sec:estimates}

Throughout this section we work on the probability space
$(\Omega, \F , \P)$ where $(\omega , \F) = (\Z^+ , 2^{\Z^+})^\infty$
and $\P$ is the probability measure making the coordinates independent
Poissons, the $n^{th}$ having mean $1/n$.  Our notation includes the
coordinate variables $\{ X_n \}$, the random multiset $M$ and its
sumset $S$.  We also define partial sums
\begin{eqnarray*}
Z_n & := & \sum_{k=1}^n X_k \, ;\\
W_n & := & \sum_{k=1}^n k X_k  \, .
\end{eqnarray*}
Thus $Z_n$ counts the cardinality of $M \cap [n]$ and $W_n$ is
the sum of all elements of $M \cap [n]$, which is the greatest
element of $\sumset (M \cap [n])$.  We will need estimates
for the right tail of $Z_n$ and $W_n$, which are obtained in a
straightforward way from their moment generating functions.

Let $\phi_{Z,n} (\lambda) := \E e^{\lambda Z_n}$ denote the moment
generating function for $Z_n$ and let $\psi_{Z,n} (\lambda)$ denote
$\log \phi_{Z,n} (\lambda)$.  Let $\tphi$ and $\tpsi$ denote the
corresponding functions for $W_n$ in place of $Z_n$.  Let
$H_n := \sum_{j=1}^n 1/j$ denote the $n^{th}$ harmonic number.
Using $\disp{\E e^{\lambda X_j} = \exp [ (e^{\lambda j} - 1) / j ]}$
and summing over $j$ leads immediately to
\begin{equation} \label{eq:psi}
\psi_{Z,n} (\lambda) = H_n \cdot \left ( e^\lambda - 1 \right ) \, .
\end{equation}
Similarly,
\begin{equation} \label{eq:tpsi}
\tpsi (\lambda) = \sum_{j=1}^n \frac{e^{j \lambda} - 1}{j}  \, .
\end{equation}
Markov's inequality implies an upper bound
\begin{equation} \label{eq:markov}
\log \P (Z_n \geq a) \leq \psi_{Z,n} (\lambda) - a \lambda
\end{equation}
for any $a > \E Z_n = H_n$.  Similarly
\begin{equation} \label{eq:markov2}
\log \P (W_n \geq a) \leq \tpsi (\lambda) - a \lambda
\end{equation}
for any $a > \E W_n = n$.

\begin{lem} \label{lem:too many} ~~\\[-3ex]
\begin{enumerate}[(i)]
\item There is a function $\beta (\ee) \sim \ee^2 / 2$ as
$\ee \downarrow 0$ such that
$$\P (Z_n \geq (1 + \ee) \log n) \leq e \, n^{- \beta (\ee)} \, .$$
\item For $\ee > 0$, let $\tau_\ee := \sup \{ n : Z_n \geq
(1 + \ee) \log n \}$.  Then $\tau_\ee < \infty$ almost surely.
\end{enumerate}
\end{lem}
\noindent{\sc Proof of Lemma \ref{lem:too many}}
For a one-sided bound one does not need to optimize~\eqref{eq:markov}
in $\lambda$ but may take the near optimal $\lambda = \log (1+\ee)$.
Set $a = (1 + \ee) \log n$ to obtain
\begin{equation} \label{eq:log eps}
\log \P (Z_n \geq (1+\ee) \log n) \;\; \leq \;\;
   H_n \ee - (1+\ee) \log (1+\ee) \log n \, .
\end{equation}
Letting $\beta (\ee) := (1 + \ee) \log (1+\ee) - \ee \sim \ee^2 / 2$
and observing that $\sup_j H_j - \log j = 1$ gives
$$\log \P (Z_n \geq (1+\ee) n ) \leq - \frac{\ee^2}{2} \log n
   + O(\ee + \ee^3 \log n)$$
which proves~$(i)$.

For~$(ii)$, apply~$(i)$ with $e^n$ in place of $n$ for $n = 1, 2, 3, \ldots$
to see that
$$\P (Z_{e^n} \geq (1 + \ee / 3) n) \leq \exp (1 - \beta (\ee / 3) n) \, .$$
By Borel-Cantelli, $Z_{e^n} \geq (1 + \ee / 3) n$ finitely often
almost surely.  For $e^{n-1} < k < e^n$, the inequality
$$\frac{Z_k}{\log k} \leq \frac{Z_{e^n}}{n-1} \leq
   \frac{n}{n-1} \frac{Z_{e^n}}{n}$$
implies that $Z_k \leq (1+\ee) \log k$ as long as $Z_{e^n}{n} \leq
1 + \ee/3$ and $n/(n-1) < 1 + \ee/3$.  We have seen by Borel-Cantelli
that these are both true for $n$ sufficiently large, proving~$(ii)$.
$\Cox$
%

The upper tail of $W_n$ may be estimated in a similar way.
Throughout the paper from this point on we will use the notation
\begin{equation} \label{eq:m}
m(n) := \lfloor n / \log n \rfloor \, .
\end{equation}

\begin{lem} \label{lem:jump} ~~\\[-3ex]
\begin{equation} \label{eq:jump}
\log \P (W_{m(n)} \geq n) \leq - \log n (\log\log n - 1) \, .
\end{equation}
It follows by Borel-Cantelli that $\tau := \sup \{ n :
W_{m(n)} \geq n \}$ is almost surely finite.
\end{lem}

\noindent{\sc Proof of Lemma \ref{lem:jump}}
The near optimal choice of $\lambda$ in~\eqref{eq:markov2} is a 
little more complicated than was the optimal choice in~\eqref{eq:markov}.
We take $\lambda := \log n \log\log n / n$ and find that
$$\log \P (W_{m(n)} \geq n) \leq \underbrace{\sum_{j=1}^{n / \log n}
   \frac{\exp(j \log n \log\log n / n) - 1}{j} }_{S_n}- \log n \log\log n \, .$$
A glance at its power series shows the function $(e^{\beta x} - 1) / x$
to be increasing in $x$ for positive $\beta$.  Hence the sum $S_n$ above may be
bounded above if we replace each term with the last term.
The number of terms is $n / \log n$ so this yields
$$\log \P (W_{m(n)} \geq n) \leq \frac{n}{\log n} \left [
   \frac{\log n - 1}{n / \log n} \right ] - \log n \log \log n
   \;\; = \;\; \log n - 1 - \log n \log \log n \, .$$
Thus $\P (W_{m(n)} \geq n) = O(n^{-\alpha})$ for any $\alpha$,
and in particular is summable.  This proves~\eqref{eq:jump}
and the summability  of $n^{1 - \log\log n}$ (since $1-\log \log n \ll 2$) finishes the
Borel-Cantelli argument.
$\Cox$



The estimates~\eqref{eq:markov}--\eqref{eq:markov2} are sharp in 
the limit when optimized over $\lambda$.  
%
\begin{lem} \label{lem:sharp}
For fixed $x > 1$, as $n \to \infty$,
$$\frac{1}{\log n} \log \P (Z_n \geq x \log n) =
   x - 1 - x \log x + o(1) \, .$$
\end{lem}

\noindent{\sc Proof of Lemma \ref{lem:sharp}:} \hfill \\
\noindent{\underline{Upper bound:}} Again we
optimize~\eqref{eq:markov}.  The optimal value of $\lambda$ is
$\log (a / H_n)$ but in fact we may use the simpler, near-optimal
value $\lambda = \log (a/\log n)$.  Setting $a = x \log n$ and
$\lambda = \log x$ yields
$$\log \P (Z_n \geq x \log n) \leq H_n (\log x - 1) - x \log n \log x$$
and plugging in $H_n = \log n + \gamma + o(1)$ yields
\begin{eqnarray*}
\log \P (Z_n \geq x \log n) & = & (\log n + \gamma + o(1)) (\log x - 1)
   - x \log x \log n \\
& = & \log n (\log x - 1 - x \log x) + (\gamma + o(1)) (\log x - 1) \, .
\end{eqnarray*}
For $1 < x \leq e$ this gives the exact upper bound
\begin{equation} \label{eq:sharp}
\log \P (Z_n \geq x \log n) \leq \log n (\log x - 1 - x \log x)
\end{equation}
while for $x > e$, one has an asymptotically negligible remainder
term on the right hand side of $(\gamma + o(1)) (\log x - 1)$.

\noindent{\underline{Lower bound:}}
We use a tilting argument.  Let $\P_x$ be the probability measure on
$(\Omega , \F)$ making the coordinates independent Poissons with means
$x/n$.  Let $G_n$ be the event that $\disp{x \log n \leq Z_n
\leq x \log n + (\log n)^{2/3}}$.  The law under $\P_x$ of $Z_n$ is
Poisson with mean $x H_n = x (\log n + O(1))$, from which it follows
that $\P_x (G_n) to 1/2$ as $n \to \infty$ for any fixed $x > 1$.
The tilting argument is simply the inequality
$$\P (G_n) \geq \P_x (G_n) \inf_{\omega \in G_n}
   \frac{d\P}{d\P_x} (\omega)$$
On the $\sigma$-field $\F_n$ generated by $X_1 , \ldots , X_n$,
the Radon-Nikodym derivative is easily computed as
\begin{eqnarray}
\frac{d\P}{d\P_x} (\omega) & = & \prod_{k=1}^n
   e^{(x - 1)/n} \; x^{-X_k} \nonumber \\[2ex]
& = & \exp ((x - 1) H_n) \; x^{-Z_n} \, . \label{eq:rhs2}
\end{eqnarray}
On the event $G_n$, we have $Z_n = x \log n + O(x)$.  Plugging
this into~\eqref{eq:rhs2} and using also $H_n = \log n + O(1)$
shows that on $G_n$,
\begin{eqnarray}
\log \frac{d\P}{d\P_x} & = & (x-1) H_n - Z_n \log x \nonumber \\[2ex]
& = & \log n (x - 1 - x \log x) + O(1) \, , \label{eq:G}
\end{eqnarray}
completing the proof.
$\Cox$


\setcounter{equation}{0}
\section{Quenched probabilities for the Poisson model and the proof
of Theorem~\protect{\ref{th:poisson}}}

The above lemmas are written for any $\ee > 0$ in case future
work requires pushing $\ee$ arbitrarily close to zero.  However,
for our purposes, $\ee = 1/100$ will be fine.  To simplify
notation (and free up $\ee$ notationally for other uses) we
set $\ee = 1/100$ in Lemma~\ref{lem:too many} and we set
$$T = \max \{ \tau_{1/100} , \tau \}$$
where $\tau$ is the supremum in Lemma~\ref{lem:jump}.
Define
\begin{eqnarray} \label{eq:pt}
p_n & := & \P \left [ n \in S \right ]  \, ; \\
\pt_n & := & \P \left [ T < m(n) \mbox{ and } n \in S \right ] \, .
\end{eqnarray}
Thus $\{ \pt_n \}$ are the so-called quenched probabilities, with
exceptional events $\{ T \geq m(n) \}$.  Although the exceptional events
vary with $n$, they form a decreasing sequence, which allows us to
assume without too much penalty that none of the exceptional
events occurs.  The following lemma encapsulates the dimension
estimate in the Poisson model.

\begin{lem}[dimension of $S$] \label{lem:summable}
There is a constant $C$ such that for all $n$,
$$\pt_n \leq C n^{-1 + \ln 2 + 0.02} \, .$$
\end{lem}

\noindent{\sc Proof of Lemma \ref{lem:summable}:}
Let $G_n$ be the event that $T < m(n)$ while also $n \in S$;
this is the event whose probability we need to bound from above.
Call a sequence $\yy = (y_1 , y_2 ,\ldots)$ admissible if it
is coordinatewise less than or equal to $\omega$.

When $G_n$ occurs, because $\tau \leq m(n)$, it is not possible
for $n$ to be a sum $\sum_j j y_j$ for an admissible $\yy$ with
$y_j$ vanishing for $j > m(n)$: even setting $y_j = \omega_j$
for $j \leq m(n)$ does not give a big enough sum.  Therefore,
breaking any admissible vector into the part below $m(n)$ and
the part at $m(n)$ or above, the event $G_n$ is contained in
the following event:
\begin{quote}
$Z_{m(n)} \leq 1.01 \log m$ and $W_{m(n)} < n$ and
there is some $k$ with $k = \sum_j j (y_j' + y_j'')$
with $\yy'$ supported on $[1,m(n)-1]$ and $y''$ nonzero
and supported on $[m(n),n]$ and both $\yy'$ and $\yy''$ admissible.
\end{quote}
Let $p_k'$ denote the probability that $Z_{m(n)} \leq 1.01 \log m$
and $W_{m(n)} \leq m$ and $k = \sum_j j y_j'$ for an admissible
$\yy'$ supported on $[1,m(n)-1]$ and let $p_k''$ be the probability
that $k = \sum_j j y_j''$ for an admissible $k$ supported on
$[m(n),n]$.  By independence of the coordinates $\omega_k$
we see we have shown that
\begin{eqnarray}
\pt_n & \leq & \sum_{k=m(n)+1}^n p_{n-k}' p_k'' \nonumber \\[1ex]
& \leq & \left ( \sum_{k=1}^{m(n)} p_k' \right ) \; \cdot \;
   \max_{m(n)+1 \leq k \leq n} p_k'' \, . \label{eq:pt bound}
\end{eqnarray}

By Fubini's theorem, the first of these factors is equal to
the expected number of $k \leq m(n)$ in $S(\omega |_{m(n)})$
where $\omega |_{m(n)}$ is the sequence $\omega$ with all
entries zeroed out above $m(n)$.  Letting $Z_m := \sum_{j=1}^{m(n)}
\omega_j$ denote the size (with multiplicity) of $\omega$ up to
$m(n)$, it is immediate that the number of $k \leq m(n)$ in
$S(\omega |_{m(n)})$ is at most $2^{Z_m}$.  But on the event
$G_n$, it always holds that $Z_m \leq 1.01 \log m$ because
$\tau_{1/100} \leq T < m(n)$.  Therefore the first factor on
the right-hand side of~\eqref{eq:pt bound} is bounded above by
\begin{equation} \label{eq:sum bound}
\sum_{k=1}^{m(n)} p_k' \leq 2^{1.01 \log m}
   = m^{1.01 \ln 2} \, .
\end{equation}

Next we claim that there is a constant $C$ such that
\begin{equation} \label{eq:max bound}
p_k'' \leq C \frac{\log^2 n}{n} \mbox{ for all } k \in [m(n),n]  \, .
\end{equation}
To prove this, start with the observation that $\P (\omega_j \geq 2)
\leq j^{-2}$ leading to
$$\P (H) \leq \sum_{j= m(n)}^n j^{-2} \leq m(n)^{-1} = \frac{\log n}{n}$$
where $H$ is the event that $\omega_j \geq 2$ for some $j \in [m(n),n]$.
The event that $k = \sum_j j y_j''$ for an admissible $k$ supported on
$[m(n),n]$ but that $H$ does not occur is contained in the union of
events $E_j$ that $\omega_j = 1$ and $k-j = \sum_i i y_i''$
for some admissible $\yy''$ supported on $[m(n),n] \setminus \{ j \}$.
Using independence of $\omega_j$ from the other coordinates of $\omega$,
along with our description of what must happen if $H$ does not, we obtain
\begin{eqnarray}
p_k'' & \leq & \P (H) + \sum_{j = m(n)}^n \frac{1}{j} p_{k-j}''
   \nonumber \\
& \leq & \P (H) + \frac{1}{m(n)} \sum_{j=m(n)}^n p_{k-j}''
   \nonumber \\
& \leq & \frac{\log n}{n} \left ( 1 + \sum_{j=m(n)}^n p_{k-j}''
   \label{eq:prod} \right) \, .
\end{eqnarray}

We may now employ a relatively easy upper bound on the summation in the
last factor, namely we may use the expected number of ways of obtaining
each number as a sum of large parts (recall that the expectation when
not restricting the parts is too large to be  useful).
Accordingly, we define the generating function
$$F (z , \omega) := \prod_t (1 + z^t)$$
where the index $t$ of the product ranges over values in
$[m(n),n]$ such that $\omega_t = 1$ (recall we ruled out
values of 2 or more).  Then
\begin{eqnarray*}
\sum_{j=m(n)}^n p_{k-j}'' & \leq & \sum_{j=m(n)}^n \E [z_j] F (z,\omega) \\
& \leq & \sum_{j=1}^\infty \E F (1,\omega) \\
& \leq & \prod_{j=m(n)}^n (1 + \frac{1}{j}) \\
& = & \frac{n+1}{m(n)} \, .
\end{eqnarray*}
Putting this together with~\eqref{eq:prod} yields
$$p_k'' \leq \frac{\log n}{n} \left ( 1 + \frac{n+1}{n} \log n \right )
   = O \left ( \frac{\log^2 n}{n} \right )$$
proving~\eqref{eq:max bound}.

Finally, plugging in~\eqref{eq:sum bound} and~\eqref{eq:max bound}
into~\eqref{eq:pt bound} shows that
$$\pt_n \leq C m(n)^{1.01 \ln 2} \;  \frac{\log^2 n}{n}$$
which is bounded about by a constant multiple of $n^{\ln 2 + 0.02 - 1}$,
completing the proof of the lemma.
$\Cox$

It is now routine to establish something that is almost
Theorem~\ref{th:poisson}.

\begin{thm} \label{th:almost}
$$\P^4 \left ( \bigcap_{r=1}^4 S(\omega^r) \mbox{ \rm is finite}
   \right ) = 1 \, .$$
\end{thm}

\noindent{\sc Proof:}
Let $T^1 , \ldots , T^4$ denote the quantities $T(\omega)$ when
$\omega = \omega_1 , \ldots , \omega_4$ respectively.
Let $T^*$ denote the maximum of $\{ T^1 , T^2 , T^3 , T^4 \}$.
By Lemma~\ref{lem:summable} and the independence of the $\omega^r$
for $1 \leq j \leq 4$, the probability that $T^r < m(n)$ and
$n \in S(\omega^r)$ for all $1 \leq r \leq 4$ is at most
a constant multiple of $n^{-3.94 + 4 \ln 2}$.  The exponent
is less than $-1$, so the series is summable and we conclude
that $n$ is in the intersection of all four sets $S(\omega^r)$
for finitely many $n > T_*$ almost surely.  Almost sure
finiteness of $T^*$ finishes the proof.
$\Cox$

\noindent{\sc Proof of Theorem}~\ref{th:poisson}:
By Theorem~\ref{th:almost} we may choose $L$ sufficiently large
so that the $\P (H_L) < 1/4$ where $E_L$ is the event
that $[L,\infty) \cap \bigcap_{r=1}^4 S(\omega^r)$
is non-empty.  The event $E_L$ is an increasing function
of the independent random variables $\{ X_j \}$.
The event $E_L'$ that $[1,L-1] \cap
\bigcap_{r=1}^4 S(\omega^r)$ is non-empty is also a
increasing function of the coordinates $X_j$ and
has some fixed nonzero probability $b_L$; for example,
$b_L$ is at least the probability that $Z_L = 0$,
which is a simple Poisson event with probability
asymptotic to $e^{-4\gamma} L^{-4}$ when $L$ is large.
Harris's inequality says that any two increasing functions
of independent random variables are nonnegatively correlated
(see, e.g.,~\cite[Section 2.2]{grimmett}).  Their complements
are also nonnegatively correlated and this gives
$$\P^4 \left [ E_L^c \cap (E_L')^c \right ] \geq
   \P(E_L^c) \P ((E_L')^c) \geq \frac{b_L}{2} \, ,$$
finishing the proof of Theorem~\ref{th:poisson}.

\setcounter{equation}{0}
\section{Computation of the marginal probabilities}

In this section we prove Theorem~\ref{th:p_k}.  This is not
needed for Theorem~\ref{th:main} but serves to establish the
so-called lottery effect, that is, the fact that $p_k$ has a
different exponent of decay from the quenched probabilities $\pt_k$.
It is generally easier to prove estimates for the partial sums
$\sum_{k=1}^n p_k$ than for $p_n$.  We take care of this first,
which is the bulk of the work.  Let $A_n := \sum_{k=1}^n p_k =
\E |S \cap [n]|$.
\begin{lem} \label{lem:A_n}
$$A_n = n^{\disp{1 + \rate + o(1)}},$$ where 
$$\rate = \frac{1 - \log 2 - \log(1/\log 2)}{\log 2}
\approx -0.08607\ldots.$$
\end{lem}

\noindent{\sc Proof:}
For the upper bound, recall the notation $Z_n = \sum_{j=1}^n X_j$
and observe that the cardinality of $|S \cap [n]|$
is at most $2^{Z_n} \wedge n$ (here $\wedge$ is used to 
denote the binary operation of taking the minimum).  This leads to
\begin{eqnarray*}
A_n & \leq & \E \left ( 2^{Z_n} \wedge n \right ) \\[1ex]
& \leq & 1 + \sum_{0 \leq k \leq \log_2 n} 2^k \P (Z_n > k)  \\[1ex]
& \leq & C \log n \sup_{0 \leq x \leq 1 / \log 2} \P (Z_n \geq x \log n)
   2^{x \log n} \, .
\end{eqnarray*}
Recalling from~\eqref{eq:sharp} that
$\log \P (Z_n \geq x \log n) \leq \log n \cdot [x - 1 - x \log x]$
we obtain
$$\frac{\log A_n}{\log n} \leq o(1) + \sup_{0 \leq x \leq 1 / \log 2} \;\;
   \left [ x \log 2 + x - 1 - x \log x \right ]  \, .$$
The supremum is achieved at the right endpoint, therefore
$$\frac{\log A_n}{\log n} \leq o(1) + \frac{1 - \log(1/\log 2)}{\log 2}
   = \rate + 1 + o(1) \, .$$

For the reverse inequality, we begin by recalling the tilted laws
$\P_x$ from the proof of the lower bound in Lemma~\ref{lem:sharp},
fixing the value $x = 1 / \log 2$ for the remainder of this proof.
The idea is that when $Z_n \approx (1 / \log 2) n$, then
$S_n \cap [n]$ should have size roughly $2^{Z_n} \approx n$.
Fix $\ee > 0$ and define
$$G_n^\ee := G_n \cap \left \{ |S_n \cap [n]| \geq n^{1-\ee}
   \right \} \, .$$
The infimum of the Radon-Nikodym derivative $d\P/d\P_x$ on $G_n$,
computed in~\eqref{eq:G}, is $n^\rate + o(1)$, so the proof is
complete once we establish
\begin{equation} \label{eq:fat}
\P_x (G_n^\ee) = 1 - o(1) \mbox{ for each fixed } \ee > 0 \, .
\end{equation}

To show~\eqref{eq:fat}, we begin with some definitions.
Let $\tau_j := \inf \{ n : Z_n = j \}$ be the $j^{th}$ smallest value
in the multiset $M$.  Let $\F_j := \sigma (X_i \wedge (j - Z_{i-1})$
be the $\sigma$-field containing the values of the $j$ elements
of the multiset $M$.  This is a natural filtration on which the
random variables $x_j := X_{\tau_j}$ form an adapted sequence.
Given $x_j$, we may easily compute
\begin{eqnarray*}
\P_x \left ( \log \frac{x_{j+1}}{x_j} > u \right ) & = &
   \left ( 1 + O(\frac{1}{x_j}) \right )  \prod_{x_j < k < x_j e^u}
   e^{- 1 / (k \log 2)} \\[1ex]
& = &  \left ( 1 + O(\frac{1}{x_j}) \right ) e^{- u / \log 2} \, ;
\end{eqnarray*}
It is not hard to see from this that the conditional distribution of 
$\log (x_{j+1}/x_j)$ given $x_j$ is stochastically bounded between 
exponentials of means $\log 2 + O(1/x_j)$, where the fudge term 
accounts for the possibility that $x+1 = x_j$ and for the discretization.
Define
\begin{eqnarray*}
s_j & := & \sumset (x_1 , \ldots , x_j) \, ; \\
Y_j & := & \log |s_j| - \log x_j \, ; \\
\Delta_j & := & Y_{j+1} - Y_j = U_j - V_j \, , \\
\mbox{where} \hspace{0.5in} && \\
U_j & := & \log |s_{j+1}| - \log |s_j| \, ; \\
V_j & := & \log x_{j+1} - \log x_j \, .
\end{eqnarray*}

\begin{lem} \label{lem:Y}
\begin{equation}
\P_x (Y_j \leq - j/4) \to 0 \mbox{ as } j \to \infty . \label{eq:prob}
\end{equation}
\end{lem}

Assuming the lemma and plugging $j = \log_2 n - \log \log n$ 
into~\eqref{eq:prob}, it follows that
$$\P_x (Y_{\log_2 n - \log \log n} \leq - \ee \log n) = o(1)$$
as $n \to \infty$ for any $\ee > 0$.  Another event whose probabilities
goes to zero is the event that $W_j \geq n$ (recall that $W_{\tau_j}$
is the sum of the elements of $M$ up to $X_{\tau_j} = x_j$).  On the
complement of this event, $s_j \subseteq S \cap [n]$.  Finally, the
event $\log_2 x_j < j - j^{2/3}$ also goes to zero.  On the
complement of the union of these three small events, $|S \cap [n]|
\geq |s_j| = x_j e^{Y_j} \geq 2^{j - j^{2/3}} n^{-\ee} \geq
n^{1 - (\log n)^{-1/3} - \ee}$.  Because $\ee > 0$ is arbitrary,
this proves~\eqref{eq:fat} and finishes the proof of Lemma~\ref{lem:A_n}
modulo Lemma~\ref{lem:Y}.
$\Cox$

The proof of Lemma~\ref{lem:Y} requires the following standard
deviation estmiate for supermartingales with bounded exponential
moment.
\begin{lem} \label{lem:tail}
Let $\{S_i\}_{i \ge 0}$ be a supermartingale with respect to the 
filtration $\{\F_i\}$, with $S_0=0$.  Suppose that the increments 
$\xi_{i+1}:=S_{i+1}-S_i$ satisfy  $\E(e^{\xi_{i+1}}| \F_i) \leq B$ 
for all $i \geq 0$.
Then for all integer $\ell>0$ and real $R \in [0, 2\ell B]$, we have
\begin{equation} \label{eq:supermart}
\P(S_\ell >R) \leq e^{-R^2/(4\ell B)} \,.
\end{equation}
\end{lem}

\noindent{\sc Proof:}
By Lemma 3.1 from \cite{freedman75},
the positive function $g(t)=(e^t-1-t)/t^2$ (where $g(0)=1/2$)  
is increasing in $\R$.  Thus  for all $\lambda \in [0,1]$ and 
$\xi \in \R$, we have
\begin{equation} \label{eq:freed}
(\lambda\xi)^2 g(\lambda\xi) \le (\lambda\xi)^2 \max\{g(0),g(\xi)\}    
   \leq \lambda^2 e^{|\xi|} \,.
\end{equation}
Because $\{S_i\}$ is a supermartingale,$\E_\ell(\xi) \le 0$ where 
(just for this proof) we abbreviate $\E_\ell(\cdot)= \E(\cdot | \F_\ell)$.  
Taking expectations in~(\ref{eq:freed}), we infer that  
$$\E_\ell(e^{\lambda\xi}) \leq 1+\lambda\E_\ell(\xi_{\ell+1})    
   + \lambda^2\E_\ell(e^{|\xi_{\ell+1}|}) \, .$$
Thus $\E_\ell (e^{\lambda \xi_{\ell+1}}) \leq 1+B\lambda^2 <
e^{B\lambda^2}$, whence $\E_\ell (e^{\lambda S_{\ell+1}}) 
\leq e^{\lambda S_{\ell}+B\lambda^2}$.  A simple induction then leads 
to $\E (e^{\lambda S_{\ell}}) \le e^{\ell B\lambda^2}$. 
We conclude that $\P(S_\ell \ge R) \le e^{\ell B\lambda^2-\lambda R}$. 
To minimize the right-hand side, we take $\lambda=R/(2\ell B)$, 
which yields the assertion of the lemma.
$\Cox$

\noindent{\sc Proof of Lemma}~\ref{lem:Y}:
The quantity $\Delta_j$ is the difference of positive variables 
$U_j$ and $V_j$.  Conditional on $\F_j$, the variable $U_j$ is
stochastically greater than $- E_j$ where $E_j$ is an exponential
of mean $\log 2 + 0.01$.  The variable $V_j$ is of necessity in the
interval $[0,\log 2]$.  We begin by showing that 
\begin{equation}
\E V_j \geq \log 2 \left ( 1 - e^{Y_j \wedge 0} \right ) \, .
\label{eq:incr} 
\end{equation}
Let $R$ denote the size of the overlap $R := |s_j \cap (x_{j+1} \oplus s_j)|$ 
where the $\oplus$ symbol in this case denotes translation of the set 
$s_j$ by $x_{j+1}$.  We may then express
$$\log \frac{s_{j+1}}{s_j} = \log \frac{2 |s_j| - R}{|s_j|}
   = \log 2 + \log \left ( 1 - \frac{R}{2 |s_j|} \right ) \, .$$
Using the fact that $R / (2 |s_j|) \in [0,1/2]$ and the bound
$\log (1 - u) \geq - u \log 4$ for $u \in [0,1/2]$ then gives
$$\E_x (V_j | \F_j) \geq \log 2 - \log 4 \frac{\E_x (R | \F_j)}{2 |s_j|}
   + O \left ( \frac{1}{x_j} \right ) \, .$$
But
$$\E_x (R | \F_j) = \sum_{a,b \in s_j} \P_x (x_{j+1} = b-a | \F_j)
   \leq |s_j|^2 / x_j$$
because $\P_x (x_{j+1} = k | \F_j) \leq 1/x_j$ for
any $k$.  Also trivially $R \leq |s_j|$, whence  
$$\frac{\E_x (R | \F_j)}{2 |s_j|} \leq \frac{s_j}{2 x_j}
   \wedge \frac{1}{2} \, .$$
Replacing $|s_j| / x_j$ by $\exp (Y_j)$ then proves~\eqref{eq:incr}.

The event $G := \{ Y_j \leq - j/4 \}$ can be covered by the
union over $0 \leq i \leq j$ of the events $G_i$ defined as follows.
Let $G_0$ be the event that for some $i \leq j$ we have
$U_i \leq - \ee j / 4$.  For $1 \leq i \leq j$ define $G_i$ 
to be the event that $Y_i \in [-\ee j/2 , -\ee j/4]$, $Y_j 
\leq - \ee j$, and $Y_t \leq - \ee/4$ for every $t \in [i,j]$.
To see that $G \subseteq \bigcup_{i=0}^j G_i$, observe that
if no jump is less than $-\ee/4$ then the last time $i \leq j$
that $Y_i \geq - \ee j /2$, we must have $Y_i \leq - \ee j /4$.

 From the fact that $U_j \geq - E_j$ and $V_j \geq 0$ we see easily that 
$$\P_x (G_0) \leq j \exp \left ( - \frac{\ee j}{4 (\log 2 + 0.01)} 
		\right ) \, .$$
A sufficient condition to imply the lemma is that there is some 
$c > 0$ such that $\P_x (G_i) < e^{-cj}$ for all $i,j$ with $1 \leq i < j$.   
This follows from an application of Lemma~\ref{lem:tail}

Fix $i$ and $j$ and for $i \leq k \leq j$ let 
$M_k = - (Y_{k \wedge \tau} - Y_i) - {k \wedge \tau} \ee / 4$, 
where $\tau$ is the least $t$ for which $Y_y \geq -\ee j/4$.  
On the event $G_i$, the value of $M_j - M_i$ is at least $R := \ee j/4$.
The expected increment $\Delta M_k := \E (M_{k+1} - M_k  \| \F_k)$ is 
zero when $k \geq \tau$ and otherwise is at most 
$\log 2 + 0.01 - \E V_k - \ee / 4$.  By~\eqref{eq:incr}, 
this is at most $0.01 + (\log 2) e^{-\ee j/4} - \ee/4$ 
which is less that zero, hence $\{ M_k \}$ is supermartingale.  
For any $\lambda < (\log 2 + 0.01)^{-1}$, and in particular for 
$\lambda = 1$, the quantity $\E e^{\lambda \Delta M_k}$ is bounded above
by some constant, $B$, independent of $i$ and $j$.  Applying 
Lemma~\ref{lem:tail} to the supermartingale $\{ S_t := M_{i+t} - M_i \}$ 
with $R = \ee j / 4$ and $\ell = j-i \leq j$, we see that
$$\P_x \left ( (M_j - M_i \geq \frac{\ee j}{4} \right ) \leq 
   \exp \left ( - \frac{(\ee j / 4)^2}{4 j B} \right )
   = \exp \left ( \frac{\ee^2}{64 B} \; j \right ) \, .$$
This completes the proof of Lemma~\ref{lem:Y}.
$\Cox$

\noindent{\it Proof of Theorem}~\ref{th:p_k}:
Typically one requires some kind of regularity to
get from an estimate on the partial sums to an estimate on the
individual summands.  Here instead we copy the proof of
Lemma~\ref{lem:summable}, using the large-index summands,
rather than some kind of monotonicity, to do the smoothing.

Recall that $m(n) := \lfloor n / \log n \rfloor$ and bound $p_n$
from below  by one minus the probability that all attempts to make
$n$ using a part of size between $n-m$ and $n$ fail.
\begin{eqnarray*}
p_n & \geq & 1 - \E \prod_{k=n-m}^n
   \left ( 1 - \frac{1}{k} \one_{n-k \in S} \right ) \\[1ex]
& \geq & 1 - \E \exp \left ( - \sum_{k=n-m}^n \frac{\one_{n-k \in S}}{k}
   \right ) \\[1ex]
& \geq & 1 - \E \exp \left ( - \frac{1}{n}
   \left | S \cap [n] \right | \right ) \, .
\end{eqnarray*}
By convexity of the exponential, the maximum value of
$\E e^{-Y/n}$ over all variables with mean $A_m$ taking values
in $[0,m]$ is achieved when $Y$ is equal to $m$ times a Bernoulli
with mean $(A_m/m)$.  This yields
\begin{equation} \label{eq:01}
p_n \geq \frac{A_m}{m} \left ( e^{m/n} - 1 \right ) \geq
   \frac{A_m}{n} \, .
\end{equation}
On the other hand, by Lemma~\ref{lem:jump},
\begin{equation} \label{eq:02}
p_n \leq n^{1 - \log \log n} + \sum_{k=m}^n \frac{1}{k} p_k
\leq \frac{\log n}{n} A_n \, .
\end{equation}
Together,~\eqref{eq:01} and~\eqref{eq:02} show that $A_n =
n^{1 + \rate + o(1)}$ implies $p_n = n^{\rate + o(1)}$;
this proves Theorem~\ref{th:p_k} modulo Lemma~\ref{lem:A_n}.
$\Cox$

\setcounter{equation}{0}
\section{Random permutations}

In this section we prove Theorem~\ref{th:main}.  The starting point
is a coupling between the Poisson variables in the Poisson model
and permutations in the group theoretic model.  The underlying space
for the coupling will be the space $(\Omega , \F , \Q)$ and its
fourfold product where
$$\Omega := \left ( \prod_{N=1}^\infty \S_N \right ) \times
   \prod_{j=1}^\infty \Z^+$$
and $\F$ is the product of the complete $\sigma$-fields
(the power set of $\S_N$ or $\Z^+$) in each coordinate.
For $\omega = (s_1 , s_2 , \ldots , x_1, x_2, \ldots)
\in \Omega$, define the coordinate functions
$X_j (\omega) := x_j' \in \Z^+$ and $\sigma_N (\omega) :=
s_N \in \S_N$.

Let $\Delta_N := ||Q_{N,m(N)} - \nu_{m(n)}||_{TV}$ be as in
the statement of Lemma~\ref{lem:AT}.  To re-iterate, $\Delta_N$
denotes the total variation distance between the joint distribution
of number of cycles of sizes $1, \ldots , m(N)$ in a uniform random
permutation from $\S_N$ and the product Poisson measure on
$(\Z^+)^{m(N)}$ whose $j^{th}$ coordinate has mean $1/j$.
The lemma of Arratia and Tavar\'e states that
\begin{equation} \label{eq:Delta}
\Delta_N \leq \exp (- C \log N \log \log N) = N^{-c \log \log N} \, .
\end{equation}

\begin{lem} \label{lem:coupling}
There is a probability measure $\Q$ on $\Omega$ such that the laws
of the random variables $\sigma_N$ and $X_j$ have the following
properties for all $N$ and $j$:
\begin{enumerate}[(i)]
\item $\sigma_N \sim \P_N$ (the uniform measure on $\S_N$);
\item $X_j \sim \pois (1/j)$ (a Poisson with mean $1/j$);
\item with probability $1 - \Delta_N$, for all $N$ and all
$j \leq m(N)$, the permutation $\sigma_N$ has exactly
$X_j$ cycles of length $j$.
\end{enumerate}
\end{lem}

\noindent{\sc Proof:}
There is a coupling $\Q_N$ of $\P_N$ and $\nu_{m(N)}$ giving measure
$1 - \Delta_N$ to the set of $(\sigma_N , \{ X_n \, : n \geq 1 \})$
such that there are $X_j$ cycles of $\sigma_N$ of length $j$
for all $j$.  The grand coupling $\Q$ may be constructed by
first making $(x_1 , x_2, \ldots)$ independent Poissons with
means $1/j$ and then giving $\sigma_N$ the conditional distribution
of $\Q_N$ given $(x_1, x_2, \ldots)$.
$\Cox$

\noindent{\sc Proof of Theorem}~\ref{th:main}:
Fix $L$ and $b_L$ as in the end of the proof of Theorem~\ref{th:poisson}.
Choose $N_0$ such that $N_0 / (\log N_0)^2 > L$.  Let
$(\Omega , \F , \Q)^4$ be the fourfold product of the measure
constructed in Lemma~\ref{lem:coupling}.  The notation is
a bit unwieldy but we will denote the generic element $\omega \in \Omega^4$
by $\langle s_j^r , x_j^r : j \in \Z^+, 1 \leq r \leq 4 \rangle$.
Let $X_j^r, 1 \leq r \leq 4$ denote the $(j,r)$ coordinate $x_j^r$
and $\sigma_j^r$ the $(j,r)$ permutation coordinate $s_j^r$ of $\omega$.
We let $\X^r$ denote the sequence $(X_j^r : j \geq 1)$.

Let $G \in \F^4$ be the event that $\bigcap_{r=1}^4 S(\X)$ is empty.
By Theorem~\ref{th:poisson} and the identification of the constant
$b_L$, we know that $\P^4 (G) \geq b_L / 2$.  Choose $N_1 \geq N_0$
so that $\Delta_{N_1} \leq b_L / 40$ and also $\P_N (T \geq
N / (\log N)^2) \leq b_L / 40$ for $N \geq N_1$.  Let $H_N^r$
denote the uncoupling event, namely the event that for some
$j \leq m(N)$, the permutation $\sigma_N^r$ has a number of
$j$-cycles different from $X_j^r$.  This has probability $\Delta_N$,
hence for $N \geq N_1$, at most $b_L / 40$.  Therefore, the event
$G_N := G \setminus (H_N^1 \cup H_N^2 \cup H_N^3 \cup H_N^4)$ has
probability at least $(2/5) b_L$.  On $G_N$, the common intersection
of $S(\sigma_N^r)$ for $a \leq r \leq 4$ cannot contain any
elements less than $m(N)$ because the cycle counts of $\sigma_N^r$
agree with $\{ X_j^r \}$ for all cycles of length at most $m(N)$
and on $G$, the sumsets of these counts have no common intersection.

\begin{lem} \label{lem:subtract}
There is an $N_2$ such that for all $N \geq N_2$,
\begin{equation} \label{eq:E}
\P^4 (E_N) \leq \frac{3 b_L}{10}
\end{equation}
where the event $E_N$ is defined by
$$E_N := G_N \cap \left \{ \bigcap_{r=1}^4 S(\sigma_N^r) \cap [m(N),N]
   \neq \emptyset \right \} \, .$$
\end{lem}

Theorem~\ref{th:main} follows from this: for any $N \geq N_2$,
$\P (G_N \setminus E_N) \geq (3/10) b_L - (1/5) b_L$.  On this
event, the four sets $S(\sigma_N^r)$ have no common intersection.
Letting $b$ be the minimum of $b_L / 10$ and the least
probability of no common intersection over all $N < N_2$
then proves the theorem.  It remains to prove the lemma.

\noindent{\sc Proof of Lemma}~\ref{lem:subtract}:
The outline is very similar to the outline of the proof of
Theorem~\ref{th:poisson}.  Fix $N \geq N_1$.  The analogue to
Lemma~\ref{lem:summable} is to define, for $m(N) \leq n \leq N$,
a quantity $\qt_n$ analogous to $\pt_n$.  This is the probability
that $n \in S(\sigma_n^r)$ while also $T^r < m(n)$; this probability
clearly does not depend on $r$.  We will show that $\qt_n^4$ is
summable.  To see that this is enough, assume it is true and pick
pick $N_2$ to make the tail sum sufficiently small:
$$\sum_{n=m(N_2)}^\infty \qt_n \leq \frac{b_L}{10} \, .$$

If $E_n$ occurs then either some $T^r \geq m(n)$ or $E_n$
occurs without this.  The first of these two probabilities
is limited to $b_L / 10$ by choice of $N_1$: $m(n)$ can be
no less than $N_1 / (\log N_1)^2$, guaranteeing that
$T^r \geq m(n)$ with probability at most $b_L / 40$
and hence $T > m(n)$ with probability at most $b_L / 10$.
The second of the two probabilities is limited to $b_L / 10$
as long as $N \geq N_2$ because the sum of $\qt_n^4$
as $n$ ranges over $[m(N),N]$ will be at most the tail sum
of $\qt_n$ from $m(N_2)$.  This makes~\eqref{eq:E} the sum of
two quantities each at most $b_L / 10$ and finishes the
proof of the lemma with one tenth to spare.

Fourth power summability of $\qt_n$ is proved via an estimate
very similar to the estimate in Lemma~\ref{lem:summable}.
Because $E_N \subseteq G_N$, the coupling is unbroken and
it is therefore not possible to have $n \in S(\sigma_n^r)$
equal to $\sum j y_j$ with $\yy$ supported on $[1,m(n)]$.
Hence, as before, $\yy$ decomposes into $\yy' + \yy''$
with $\yy'$ supported on $[1,m(n)]$ and $\yy'$ supported
on $[m(n)+1,n]$ and not identically zero.  Also as before
we have the upper bound
\begin{equation} \label{eq:qt bound}
\qt_n \leq \left ( \sum_{k=1}^{m(n)} q_k' \right ) \; \cdot \;
   \max_{m(n)+1 \leq k \leq n} q_k''
\end{equation}
where $q_k'$ is the probability $\qt_k$ but using only
cycles of size at most $m(n)$ and $q_k''$ is the analogue
of $\qt_k$ when only cycles of size at least $m(n)+1$ are used.

Analogously to~\eqref{eq:sum bound}, the first factor is at most
$m(n)^{1.01 \log 2}$ because the coupling is unbroken and we already
proved this bound for the Poisson variables.  It suffices therefore to
prove the bound
\begin{equation} \label{eq:max q}
q_n'' \leq C \frac{\log^3 n}{n} \mbox{ for } n \geq m(N)
\end{equation}
analogous to~\eqref{eq:max bound}.  Here the proof diverges from
the proof of Lemma~\ref{lem:summable} because the constraint
on the vector $\yy''$ in $\sum j y_j''$ is that $y_j''$ be at most
the number $Y_j$ of $j$-cycles in the actual permutation $\sigma_N^r$,
rather than being at most the Poisson variable $X_j^r$.  The variables
$Y_j$ as were the variables $X_j^r$ are not independent so instead
we argue as follows.

Recall that $N_2 / (\log N_2)^2 \leq N / (\log N)^2 \leq m(n) \leq N$
and observe that the quantity $q_n''$ is at most the sum over
$j \geq m(n)$ of
$$\P_N \left [ Y_j \geq 1 \mbox{ and } n-j = \sum_i i y_i''
   \mbox{ for some } \yy'' \leq \Y - \delta_j
   \mbox{ supported on } [m(n)+1,N] \right ] \, .$$
Here we have denoted by $\Y$ the vector whose components are $Y_j$.
The actual elements in the cycles of the permutation $\sigma_N^r$ are
exchangeable given the cycle lengths, so the probability that
the element~1 is in the cycle of length $j$ is at least
$N / m(n)$, which is a least $(\log n)^{-2}$.  Therefore,
we may write a new upper bound
$$q_n'' \leq (\log n)^2 \sum_{j \geq m(n)} \P_N \left [ 1 \mbox{ is
   in a cycle of length } j \mbox{ and } n-j \in S^* (\sigma_N^r)
   \vphantom{Y^{B^i}} \right ]$$
where the superscript $S^*$ denotes that we count only sums of
cycle sizes at least $m(n)$.

The reason for going through this trouble is that conditioned
on~1 being in a cycle of size $j$, the remainder of the
permutation is uncontaminated: its cycle sizes are distributed
as those of a uniform pick from $\S_{n-j}$.  Also, the probability
of~1 being in a cycle of size $j$ is precisely $1/N$.  Therefore,
\begin{equation} \label{eq:star}
q_n'' \leq \frac{(\log n)^2}{N} \sum_{j \geq m(n)}
   \P_N \left [ n-j \in S^* P(\sigma_{N-j}^r)
   \vphantom{Y^{B^i}} \right ] \, .
\end{equation}

To evaluate the summand, first observe that replacing $S^*$ with $S$,
the expected number of invariant sets of $\sigma_{N-j}^r$ of size
$n-j$ is precisely~1 for any $N, n, j$.  Next, consider any
invariant set of $\sigma_{N-j}^r$ of size $n-j$ and bound
from above the probability that it consists entirely of
cycles larger than $m(n)$.  This is the same as the probability
that a random element of $\S_{n-j}$ has only cycles of
length at least $m(n)$.  We may evaluate this via the
Arratia-Tavar\'e lemma: it is at most equal to $\P (Z_{m(n)} = 0)$
(the probability that a Poisson ensemble with $\E X_j = 1/j$
takes only value zero up to $j = m(n)$) plus the total
variation distance between the Poisson product measure
and the actual counts of cycle sizes up to $m(n)$.
By the bound in Lemma~\ref{lem:AT}, this total variation
distance is at most $\exp [C  n/m(n) \log (n/m(n)) ]$
which is bounded above by $n^{-C \log \log n}$ and hence
decays faster than any polynomial in $n$.  The probability
of $Z_{m(n)} = 0$ is $e^{-H_{m(n)}} \leq 1/m(n)$, whence
the summand in~\eqref{eq:star} is therefore at most
$1/m(n) \sim \log n / n$ and the whole sum is at most
$\log n (N/n)$.  This makes the right-hand side of~\eqref{eq:star}
at most $\log^3 n / n$, establishing~\eqref{eq:max q} and
completing the proof of Lemma~\ref{lem:subtract} and hence
of Theorem~\ref{th:main}.
$\Cox$
\bibliographystyle{alpha}
\bibliography{RPP}
\end{document}